\documentclass[oneside,french]{amsart}
\usepackage[T1]{fontenc}
\usepackage[latin9]{inputenc}
\usepackage{amsbsy}
\usepackage{amstext}
\usepackage{amsthm}
\usepackage{amssymb}
\usepackage[all]{xy}

\makeatletter

\providecommand{\tabularnewline}{\\}

\numberwithin{equation}{section}
\numberwithin{figure}{section}
\theoremstyle{plain}
\newtheorem{thm}{\protect\theoremname}
  \theoremstyle{plain}
  \newtheorem{prop}[thm]{\protect\propositionname}
  \theoremstyle{remark}
  \newtheorem{rem}[thm]{\protect\remarkname}

\newcommand{\xyR}[1]{\xydef@\xymatrixrowsep@{#1}}
\newcommand{\xyC}[1]{\xydef@\xymatrixcolsep@{#1}}

\makeatother

\usepackage{babel}
\makeatletter
\addto\extrasfrench{%
   \providecommand{\fg}{\ifdim\lastskip>\z@\unskip\fi~\frqq}%
}

\makeatother
  \providecommand{\propositionname}{Proposition}
  \providecommand{\remarkname}{Remarque}
\providecommand{\theoremname}{Théorème}

\begin{document}

\title{Sur les triangles avec un côté minuscule}
\begin{abstract}
Soit $W\subset O(V)$ le groupe de Weyl d'un système de racines $R\subset V$.
Si $a+b+c=0=a'+b'+c'$ avec $a$, $b$ et $c$ respectivement conjugués
à $a'$, $b'$ et $c'$ dans $V$, alors $(a,b,c)$ est conjugué à
$(a',b',c')$ dans $V^{3}$ lorsque $a$ est, dans chaque composante
irréductible de $V$, colinéaire à un copoids minuscule. 
\end{abstract}


\author{Christophe Cornut}
\maketitle

\section{Introduction}

Dans un espace affine euclidien $A$, deux triangles sont transformés
l'un dans l'autre par une isométrie de $A$ si et seulement si les
longueurs de leurs côtés sont les mêmes, c'est-à-dire si et seulement
si leurs côtés sont, deux à deux, transformés l'un dans l'autre par
une isométrie de $A$ \cite[Proposition 8]{Eu90}. Pour l'espace vectoriel
euclidien sous-jacent $V$ de $A$, cette propriété élémentaire se
traduit ainsi: l'injection 
\[
V_{0}^{3}=\left\{ (a,b,c)\in V^{3}:a+b+c=0\right\} \hookrightarrow V^{3}
\]
induit une \emph{injection} de $G\backslash V_{0}^{3}$ dans $(G\backslash V)^{3}$,
où $G$ est le groupe orthogonal $O(V)$. Cela n'est plus nécessairement
vrai pour un sous-groupe $G$ de $O(V)$, comme on voit en dimension
$2$ avec $G=SO(V)$. Lorsque $G=W$ est le groupe de Weyl d'un système
de racines réduit $R$ dans $V$, la propriété ci-dessus reste cependant
vraie pour les triangles dont l'un des côtés est \emph{minuscule};
on dit d'un élément $a$ de $V$ qu'il est minuscule lorsque $\left|\{(\alpha,a):\alpha\in R_{i}\}\right|\leq3$
pour toute composante irréductible $R_{i}$ de $R$, où $(-,-)$ est
le produit scalaire $W$-invariant de $V$. Ainsi:
\begin{prop}
Si $a\in V$ est minuscule, alors pour toute paire $(b,b')$ d'éléments
$W$-conjugués de $V$ telle que $a+b$ et $a+b'$ sont également
$W$-conjugués, $b$ et $b'$ sont déjà conjugués sous le stabilisateur
$W_{a}$ de $a$ dans $W$. 
\end{prop}

\begin{rem}
Si $a$ n'est pas minuscule, le résultat peut être vrai ou faux. Pour
le système de racines de type $A_{2}$ dans $V=\{(x,y,z)\in\mathbb{R}^{3}:x+y+z=0\}$,
où $W=\mathfrak{S}_{3}$, le résultat est faux pour $a=(-1,0,1)$:
$b=(1,-1,0)$ et $b'=(0,1,-1)$ sont $W$-conjugués, $a+b=(0,-1,1)$
et $a+b'=(-1,1,0)$ sont $W$-conjugués, mais $b$ et $b'$ ne sont
pas conjugués sous le stabilisateur $W_{a}=\{1\}$ de $a$ dans $W$.
Inversement, pour les systèmes de racines de type $B_{n}$ ou $C_{n}$
dans $V=\mathbb{R}^{n}$, qui ont le même groupe de Weyl mais pas
les mêmes éléments minuscules, le résultat reste évidemment vrai pour
tout $a\in V$ qui est minuscule pour le système de racines dual,
$C_{n}$ ou $B_{n}$. 
\end{rem}

\begin{prop}
Si $a\in V$ est minuscule, alors l'algèbre des fonctions polynomiales
$W_{a}$-invariantes sur $V$ est engendrée par les fonctions polynomiales
$W$-invariantes $f:V\rightarrow\mathbb{R}$ et leurs translatées
$v\mapsto f(v+a)$. 
\end{prop}

Pour tout groupe fini $G$ agissant linéairement sur $V$, les $G$-orbites
dans $V$ sont séparées par les fonctions polynomiales $G$-invariantes
sur $V$. La proposition $1$ résulte donc de la proposition $2$.
Pour leur démonstration, on se ramène facilement au cas où $R$ est
irréductible et engendre $V$. On peut aussi supposer que $a$ est
non nul, dominant relativement au choix d'une base $\Delta$ de $R$,
et quitte à multiplier $a$ par un scalaire strictement positif, on
peut enfin supposer que $\{(\alpha,a):\alpha\in R\}\subset\{0,\pm1\}$.
On dit alors que $a$ est un \emph{copoids dominant minuscule} de
$R$. Ce sont ceux des éléments de la base duale de $\Delta$ qui
correspondent aux racines simples de multiplicité $1$ dans la plus
haute racine positive de $R$. Le tableau ci-dessous, extrait de \cite{BoLie46},
donne pour chaque système de racines irréductible réduit les multiplicités
des racines simples dans la plus haute racine positive, et les racines
qui nous concernent y ont été encadrées. On traite séparément chacun
de ces cas dans les sections suivantes, en s'appuyant sur les résultats
connus concernant les invariants de $W$ dans l'algèbre symétrique
$S$ de $V$ (\cite{BoLie46} pour $A_{n}$, $B_{n}$, $C_{n}$ et
$D_{n}$, et \cite{Me88} pour $E_{6}$ et $E_{7}$). On y identifie
systématiquement $S$ à l'algèbre des fonctions polynomiales sur $V$
grâce au produit scalaire de $V$. La translation $f(v)\mapsto f(v+a)$
sur les fonctions devient l'automorphisme $\tau(a)$ de $S$ qui envoie
$v\in V\subset S$ sur $v+(v,a)$. Il s'agit de montrer que la sous-algèbre
$\boldsymbol{S}_{a}$ de $S^{W_{a}}$ engendrée par $S^{W}$ et $\tau(a)(S^{W})$
est égale à $S^{W_{a}}$. \xyR{1pc} \xyC{1pc}
\begin{center}
\begin{tabular}{|c|c|}
\hline 
$R$ & Diagramme de Dynkin avec multiplicités\tabularnewline
\hline 
\hline 
$A_{n}$ & $\xymatrix{*+[F]{1}\ar@{-}[r] & *+[F]{1}\ar@{-}[r] & \cdots\ar@{-}[r] & *+[F]{1}}
$\vphantom{$\frac{\frac{A}{B}}{\frac{A}{B}}$}\tabularnewline
\hline 
$B_{n}$ & $\xymatrix{*+[F]{1}\ar@{-}[r] & 2\ar@{-}[r] & \cdots\ar@{-}[r] & 2\ar@{=}|-{>}[r] & 2}
$\vphantom{$\frac{\frac{A}{B}}{\frac{A}{B}}$}\tabularnewline
\hline 
$C_{n}$ & $\xymatrix{2\ar@{-}[r] & \cdots\ar@{-}[r] & 2\ar@{=}|-{<}[r] & *+[F]{1}}
$\vphantom{$\frac{\frac{A}{B}}{\frac{A}{B}}$}\tabularnewline
\hline 
$D_{n}$ & $\xymatrix{*+[F]{1}\ar@{-}[r] & 2\ar@{-}[r] & \cdots\ar@{-}[r] & 2\ar@{-}[r] & *+[F]{1}\\
 &  &  & *+[F]{1}\ar@{-}[u] & \,
}
$\vphantom{$\frac{\frac{A}{B}}{\frac{A}{B}}$}\tabularnewline
\hline 
$E_{8}$ & $\xymatrix{2\ar@{-}[r] & 3\ar@{-}[r] & 4\ar@{-}[r] & 5\ar@{-}[r] & 6\ar@{-}[r] & 4\ar@{-}[r] & 2\\
 &  &  &  & 3\ar@{-}[u]
}
$\vphantom{$\frac{\frac{A}{B}}{\frac{A}{B}}$}\tabularnewline
\hline 
$E_{7}$ & $\xymatrix{*+[F]{1}\ar@{-}[r] & 2\ar@{-}[r] & 3\ar@{-}[r] & 4\ar@{-}[r] & 3\ar@{-}[r] & 2\\
 &  &  & 2\ar@{-}[u]
}
$\vphantom{$\frac{\frac{A}{B}}{\frac{A}{B}}$}\tabularnewline
\hline 
$E_{6}$ & $\xymatrix{*+[F]{1}\ar@{-}[r] & 2\ar@{-}[r] & 3\ar@{-}[r] & 2\ar@{-}[r] & *+[F]{1}\\
 &  & 2\ar@{-}[u]
}
$\vphantom{$\frac{\frac{A}{B}}{\frac{A}{B}}$}\tabularnewline
\hline 
$F_{4}$ & $\xymatrix{2\ar@{-}[r] & 3\ar@{=}|-{>}[r] & 4\ar@{-}[r] & 2}
$\vphantom{$\frac{\frac{A}{B}}{\frac{A}{B}}$}\tabularnewline
\hline 
$G_{2}$ & $\xymatrix{2\ar@3{-}|-{>}[r] & 3}
$\vphantom{$\frac{\frac{A}{B}}{\frac{A}{B}}$}\tabularnewline
\hline 
\end{tabular}
\par\end{center}

Je suis arrivé à cette étrange question via la théorie de Hodge $p$-adique,
en étudiant une stratification des variétés affines de Deligne-Lusztig
basée sur une classification de triangles, qui me semblait plus naturelle
que les stratifications connues, basées sur la classification des
côtés de ces triangles. Les résultats de cet article montrent que
les deux approches sont équivalentes dans le cas minuscule. Je remercie
Ariane Mezard de m'avoir encouragé à traiter les cas exceptionnels. 

\section{Le cas de $A_{n}$}

Soit $n$ un entier positif. Références pour cette section: \cite[Chap. VI, \S 4, n°7]{BoLie46}.

\subsection{~}

On munit $V_{n+1}=\mathbb{R}^{n+1}$ du produit scalaire $((x_{i}),(y_{i}))=\sum x_{i}y_{i}$
et de la forme quadratique $q(x)=\frac{1}{2}(x,x)$. On note $I_{n+1}=\mathbb{Z}^{n+1}$,
$e=(1^{n+1})\in I_{n+1}$, $V_{n}$ l'orthogonal de $e$, et $A_{n}=I_{n+1}\cap V_{n}=\{(x_{i})\in I_{n+1}:\sum x_{i}=0\}$.
C'est un réseau pair d'indice $n+1$ dans son dual $A_{n}^{\ast}\subset V_{n}$,
qui est engendré par ses racines $R_{n}=\{\alpha\in A_{n}:q(\alpha)=1\}$,
de base et diagramme de Dynkin donnés par $\xyR{1pc}$$\xyC{1pc}$
\[
\begin{array}{cccccc}
\alpha_{1} & 1 & -1\\
\alpha_{2} &  & 1 & -1\\
\vdots &  &  & \ddots & \ddots\\
\alpha_{n} &  &  &  & 1 & -1
\end{array}\qquad\xymatrix{\alpha_{1}\ar@{-}[r] & \alpha_{2}\ar@{-}[r] & \cdots\ar@{-}[r] & \alpha_{n}}
\]
La plus haute racine est $\alpha_{1}+\cdots+\alpha_{n}$. Pour $r\in\{1,\cdots,n\}$
et $s=n+1-r$, 
\[
a_{r}=\frac{1}{n+1}\left(\underbrace{s,\cdots,s}_{r\,\text{termes}},\underbrace{-r,\cdots,-r}_{s\,\text{termes}}\right)\in A_{n}^{\ast}
\]
est dominant minuscule, orthogonal à toutes les racines simples sauf
$\alpha_{r}$, avec $q(a_{r})=\frac{rs}{2(n+1)}$. On note $W_{n}$
le groupe de Weyl de $R_{n}$. On étend l'action de $W_{n}$ sur $V_{n}$
à $V_{n+1}=V_{n}\bot\mathbb{R}e$ par l'action triviale sur $\mathbb{R}e$.
On obtient ainsi l'action usuelle de $W_{n}=\mathfrak{S}_{n+1}$ sur
$V_{n+1}$, par permutation des coordonnées. Le stabilisateur de $a_{r}$
dans $W_{n}$ est le sous-groupe $W_{r-1}\times W_{s-1}=\mathfrak{S}_{r}\times\mathfrak{S}_{s}$
qui permute entre elles les $r$-premières et les $s$-dernières coordonnées.
Il agit sur l'orthogonal $V_{r-1}\bot V_{s-1}$ de $a_{r}$ dans $V_{n}$,
où l'on a posé 
\begin{align*}
V_{r-1} & =\left\{ (x_{1},\cdots,x_{r},0^{s}):{\textstyle \sum}x_{i}=0\right\} ,\\
\text{et}\quad V_{s-1} & =\left\{ (0^{r},y_{1},\cdots,y_{s}):{\textstyle \sum}y_{i}=0\right\} .
\end{align*}

\subsection{~\label{subsec:GenAn}}

On note $S_{n}\subset S_{n+1}$ les algèbres symétriques de $V_{n}\subset V_{n+1}$,
donc
\[
S_{n}[e]=S_{n+1}\quad\text{et}\quad S_{n}^{W_{n}}[e]=S_{n+1}^{\mathfrak{S}_{n+1}}.
\]
Il est bien connu que $S_{n+1}^{\mathfrak{S}_{n+1}}$ est une algèbre
graduée de polynômes engendrée par les polynômes symétriques élémentaires
en $v\in\mathcal{V}$. D'après les relations de Newton de~\cite[Chap. IV, §6, n°4, Lemme 4]{BoAl4a7}
ou \cite[§3]{Me88}, cette algèbre est aussi engendrée par 
\[
\mathfrak{v}_{i}=\sum_{v\in\mathcal{V}}v^{i}\quad\text{pour}\quad i=\{1,\cdots,n+1\}
\]
où $\mathcal{V}$ est la $\mathfrak{S}_{n+1}$-orbite de $v_{1}=(1,0,\cdots,0)$
dans $V_{n+1}$. Comme $e=\mathfrak{v}_{1}$ est fixé par $W_{n}=\mathfrak{S}_{n+1}$
et $a_{1}=v_{1}-\frac{1}{n+1}e$ est la projection orthogonale de
$v_{1}$ sur $V_{n}$, $\mathcal{V}=A+\frac{1}{n+1}e$ où $A$ est
la $W_{n}$-orbite de $a_{1}$ dans $V_{n}$, donc
\[
\mathfrak{v}_{i}=\sum_{a\in A}(a+{\textstyle \frac{1}{n+1}}e)^{i}=\sum_{j=0}^{i}\left(\begin{smallmatrix}i\\
j
\end{smallmatrix}\right)\mathfrak{a}_{j}\cdot({\textstyle \frac{1}{n+1}}e)^{i-j}\quad\text{dans}\quad S_{n}^{W_{n}}[e]
\]
avec $\mathfrak{a}_{i}=\sum_{a\in A}a^{i}$ dans $S_{n}$. On a donc
aussi $S_{n+1}^{\mathfrak{S}_{n+1}}=\mathbb{R}[\mathfrak{a}_{2},\cdots,\mathfrak{a}_{n+1}][e]$
et 
\[
S_{n}^{W_{n}}=\mathbb{R}[\mathfrak{a}_{2},\cdots,\mathfrak{a}_{n+1}].
\]

\subsection{~}

La $W_{n}$-orbite $A\subset V_{n}$ de $a_{1}$ se décompose en deux
$W_{r-1}\times W_{s-1}$-orbites, $B'$ et $C'$, de cardinal $r$
et $s$, avec $a_{1}\in B'$ et 
\[
B'=\left\{ a\in A:(a,a_{r})=\frac{s}{n+1}\right\} \quad\text{et}\quad C'=\left\{ a\in A:(a,a_{r})=\frac{-r}{n+1}\right\} .
\]
Les projections orthogonales $B$ et $C$ de ces orbites sur $V_{r-1}\bot V_{s-1}\subset V_{n}$
sont
\[
B=B^{\prime}-{\textstyle \frac{1}{r}}a_{r}=W_{r-1}\cdot{\textstyle \frac{1}{r}}\left(r-1,(-1)^{r-1},0^{s}\right)\subset V_{r-1}
\]
\[
C=C'+{\textstyle \frac{1}{s}}a_{r}=W_{s-1}\cdot{\textstyle \frac{1}{s}}\left(0^{r},s-1,(-1)^{s-1}\right)\subset V_{s-1}
\]
Comme ci-dessus, on a donc
\[
S_{r-1}^{W_{r-1}}=\mathbb{R}[\mathfrak{b}_{2},\cdots,\mathfrak{b}_{r}]\quad\text{et}\quad S_{s-1}^{W_{s-1}}=\mathbb{R}[\mathfrak{c}_{2},\cdots,\mathfrak{c}_{s}]
\]
où $\mathfrak{b}_{i}=\sum_{b\in B}b^{i}$ et $\mathfrak{c}_{i}=\sum_{c\in C}c^{i}$.
Comme $S_{n}=S_{r-1}\otimes S_{s-1}[a_{r}]$, on obtient 
\[
S_{n}^{W_{r-1}\times W_{s-1}}=S_{r-1}^{W_{r-1}}\otimes S_{s-1}^{W_{s-1}}[a_{r}]=\mathbb{R}[\mathfrak{b}_{2},\cdots,\mathfrak{b}_{r},\mathfrak{c}_{2},\cdots,\mathfrak{c}_{s},a_{r}].
\]

\subsection{~}

On note $\tau=\tau(a_{r})\in\mathrm{Aut}(S_{n})$ et $\boldsymbol{S}$
la sous-algèbre de $S_{n}^{W_{r-1}\times W_{s-1}}$ engendrée par
$S_{n}^{W_{n}}$ et $\tau(S_{n}^{W_{n}})$. On se propose de montrer
que $\boldsymbol{S}=S_{n}^{W_{r-1}\times W_{s-1}}$, et il suffit
pour cela de vérifier que $a_{r}$ et tous les $\mathfrak{b}_{i}$,
$\mathfrak{c}_{i}$ sont dans $\boldsymbol{S}$.

\subsection{~}

Soit $\mathfrak{r}_{2}=\sum_{\alpha\in R_{n}}\alpha^{2}$. C'est un
élément de degré $2$ de $S_{n}^{W_{n}}$ et 
\[
(\tau-1)(\mathfrak{r}_{2})=\sum_{\alpha\in R_{n}}(\alpha+(a_{r},\alpha))^{2}-\alpha^{2}=2\sum_{\alpha\in R_{n}^{+}}(\alpha+1)^{2}-\alpha^{2}=4\sum_{\alpha\in R_{n}^{+}}\alpha+2\left|R_{n}^{+}\right|
\]
où $R_{n}^{+}=\{\alpha\in R_{n}:(\alpha,a_{r})=1\}$ est de cardinal
$rs$ et stable sous $W_{r-1}\times W_{s-1}$. Puisque $V_{n}^{W_{r-1}\times W_{s-1}}=\mathbb{R}a_{r}$,
$\sum_{\alpha\in R_{n}^{+}}\alpha=\lambda a_{r}$ avec $\lambda(a_{r},a_{r})=\sum_{R_{n}^{+}}(a_{r},\alpha)=\left|R_{n}^{+}\right|$,
donc $\lambda=n+1$. Ainsi, $4(n+1)a_{r}=(\tau-1)(\mathfrak{r}_{2})-2rs$
et $a_{r}\in\boldsymbol{S}$.

\subsection{~}

Comme $A=B'\coprod C'$, $\mathfrak{a}_{i}=\mathfrak{b}_{i}^{\prime}+\mathfrak{c}_{i}^{\prime}$
avec $\mathfrak{b}_{i}^{\prime}=\sum_{b\in B'}b^{i}$ et $\mathfrak{c}_{i}^{\prime}=\sum_{c\in C'}c^{i}$.
On a 
\[
(\tau-1)\mathfrak{a}_{i+1}=\sum_{j=0}^{i}\left(\begin{smallmatrix}i+1\\
j
\end{smallmatrix}\right)\left[\left({\textstyle \frac{s}{n+1}}\right)^{i+1-j}\mathfrak{b}_{j}^{\prime}+\left({\textstyle \frac{-r}{n+1}}\right)^{i+1-j}\mathfrak{c}_{j}^{\prime}\right].
\]
Or $\mathfrak{a}_{i}$ et $(\tau-1)\mathfrak{a}_{i+1}$ sont dans
$\boldsymbol{S}$. Puisque la matrice
\[
\left(\begin{array}{cc}
1 & 1\\
{\textstyle \frac{s}{n+1}} & {\textstyle \frac{-r}{n+1}}
\end{array}\right)
\]
est inversible, on en déduit par récurrence sur $i$ que $\mathfrak{b_{i}^{\prime}}$
et $\mathfrak{c}_{i}^{\prime}$ sont dans $\boldsymbol{S}$ pour tout
$i\in\mathbb{N}$. Or $B=B^{\prime}-{\textstyle \frac{1}{r}}a_{r}$
et $C=C^{\prime}+{\textstyle \frac{1}{s}}a_{r}$, donc
\[
\mathfrak{b}_{i}=\sum_{j=0}^{i}\left(\begin{smallmatrix}i\\
j
\end{smallmatrix}\right)\left(-{\textstyle \frac{1}{r}}a_{r}\right)^{i-j}\mathfrak{b}_{j}^{\prime}\quad\text{et}\quad\mathfrak{c}_{i}=\sum_{j=0}^{i}\left(\begin{smallmatrix}i\\
j
\end{smallmatrix}\right)\left({\textstyle \frac{1}{s}}a_{r}\right)^{i-j}\mathfrak{c}_{j}^{\prime}
\]
sont également dans $\boldsymbol{S}$.

\section{Les cas $D_{n}$, $B_{n}$ et $C_{n}$}

Soit $n$ un entier positif, qui sera supérieur ou égal à $2$ dans
les cas $B_{n}$ et $C_{n}$, et supérieur ou égal à $3$ dans le
cas $D_{n}$. Références: \cite[Chap. VI, \S 4, n°5, 6, 8]{BoLie46}.

\subsection{~}

On munit $V_{n}=\mathbb{R}^{n}$ du produit scalaire $((x_{i}),(y_{i}))=\sum x_{i}y_{i}$
et de la forme quadratique $q(x)=\frac{1}{2}(x,x)$. On note $I_{n}=\mathbb{Z}^{n}$
et $D_{n}=\{(x_{i}):\sum x_{i}\equiv0\bmod2\}$. C'est un réseau pair
d'indice $4$ dans son dual $D_{n}^{\ast}\subset V_{n}$, qui est
engendré par ses racines $R_{n}^{D}=\{\alpha\in D_{n}:q(\alpha)=1\}$,
de base et diagramme de Dynkin donnés par
\[
\begin{array}{cccccc}
\alpha_{1} & 1 & -1\\
\alpha_{2} &  & 1 & -1\\
\vdots &  &  & \ddots & \ddots\\
\alpha_{n-1} &  &  &  & 1 & -1\\
\alpha_{n}^{D} &  &  &  & 1 & 1
\end{array}\qquad\xymatrix{\alpha_{1}\ar@{-}[r] & \alpha_{2}\ar@{-}[r] & \cdots\ar@{-}[r] & \alpha_{n-2}\ar@{-}[r] & \alpha_{n-1}\\
 &  &  & \alpha_{n}^{D}\ar@{-}[u]
}
\]
La plus haute racine est $\alpha_{1}+2\alpha_{2}+\cdots+2\alpha_{n-2}+\alpha_{n-1}+\alpha_{n}^{D}$
et 
\begin{align*}
b & =(1,0,\cdots,0)\\
c' & ={\textstyle \frac{1}{2}}(1,\cdots,1,-1)\\
c & ={\textstyle \frac{1}{2}}(1,\cdots,1)
\end{align*}
sont les poids dominants minuscules, orthogonaux à toutes les racines
simples sauf respectivement $\alpha_{1}$, $\alpha_{n-1}$ et $\alpha_{n}^{D}$.
On a $q(b)=\frac{1}{2}$ et $q(c)=q(c')=\frac{n}{8}$. Le groupe de
Weyl est $W_{n}^{\circ}=\{\pm1\}_{\circ}^{n}\rtimes\mathfrak{S}_{n}$,
où $\{\pm1\}_{\circ}^{n}=\{(\epsilon_{i})\in\{\pm1\}^{n}:\prod\epsilon_{i}=1\}$. 

\subsection{~}

C'est un sous-groupe distingué de $W_{n}=\{\pm1\}^{n}\rtimes\mathfrak{S}_{n}$,
qui est le groupe de Weyl commun aux deux systèmes de racines duaux
de type $B_{n}$ et $C_{n}$, 
\begin{align*}
R_{n}^{B} & =\{\pm\epsilon_{i}\pm\epsilon_{j}:1\leq i<j\leq n\}\cup\{\pm\epsilon_{i}:1\leq i\leq n\}\\
R_{n}^{C} & =\{\pm\epsilon_{i}\pm\epsilon_{j}:1\leq i<j\leq n\}\cup\{\pm2\epsilon_{i}:1\leq i\leq n\}
\end{align*}
où $\epsilon_{i}=(\delta_{i,j})_{j}$ est la base canonique. Des bases
de ces systèmes sont données par 
\[
\begin{array}{cccccc}
\alpha_{1} & 1 & -1\\
\alpha_{2} &  & 1 & -1\\
\vdots &  &  & \ddots & \ddots\\
\alpha_{n-1} &  &  &  & 1 & -1\\
\alpha_{n}^{B} &  &  &  &  & 1
\end{array}\quad\text{et}\quad\begin{array}{cccccc}
\alpha_{1} & 1 & -1\\
\alpha_{2} &  & 1 & -1\\
\vdots &  &  & \ddots & \ddots\\
\alpha_{n-1} &  &  &  & 1 & -1\\
\alpha_{n}^{C} &  &  &  &  & 2
\end{array}
\]
Leurs diagrammes de Dynkin sont$\xyC{1.5pc}$
\[
\xymatrix{\alpha_{1}\ar@{-}[r] & \cdots\ar@{-}[r] & \alpha_{n-1}\ar@{=}|-{>}[r] & \alpha_{n}^{B}}
\quad\text{et}\quad\xymatrix{\alpha_{1}\ar@{-}[r] & \cdots\ar@{-}[r] & \alpha_{n-1}\ar@{=}|-{<}[r] & \alpha_{n}^{C}}
,
\]
leurs plus hautes racines sont
\[
\alpha_{1}+2\alpha_{2}+\cdots+2\alpha_{n-1}+2\alpha_{n}^{B}\quad\text{et}\quad2\alpha_{1}+\cdots+2\alpha_{n-1}+\alpha_{n}^{C},
\]
et leurs poids dominants minuscules sont 
\[
b=(1,0,\cdots,0)\quad\text{et}\quad c={\textstyle \frac{1}{2}}(1,\cdots,1),
\]
orthogonaux à toutes les racines simples sauf, respectivement, $\alpha_{1}$
et $\alpha_{n}^{C}$. On note $W_{b}\simeq\{\pm1\}^{n-1}\rtimes\mathfrak{S}_{n-1}$
et $W_{c}\simeq\mathfrak{S}_{n}$ les stabilisateurs, $V_{b}=\{(x_{i}):x_{1}=0\}$
et $V_{c}=\{(x_{i}):\sum x_{i}=0\}$ les orthogonaux de $b$ et $c$
dans $V_{n}$, et $S_{b},S_{c}\subset S_{n}$ les algèbres symétriques
de $V_{b},V_{c}\subset V_{n}$. On note $\beta=\tau(b)$ et $\gamma=\tau(c)$,
$\boldsymbol{S}_{b}$ la sous-algèbre de $S_{n}^{W_{b}}$ engendrée
par $S_{n}^{W_{n}}$ et $\beta(S_{n}^{W_{n}})$, et $\boldsymbol{S}_{c}$
la sous-algèbre de $S_{n}^{W_{c}}$ engendrée par $S_{n}^{W_{n}}$
et $\gamma(S_{n}^{W_{n}})$. On se propose de montrer que $\boldsymbol{S}_{b}=S_{n}^{W_{b}}$
et $\boldsymbol{S}_{c}=S_{n}^{W_{c}}$.

\subsection{~}

On note $B$ la $W_{n}$-orbite de $b$, qui a $2n$ éléments. Elle
se décompose en trois $W_{b}$-orbites, les points fixes $\pm b$
et leur complémentaire $B'$, qui est orthogonal à $b$, et en deux
$W_{c}$-orbites $B^{+}$ et $B^{-}=-B^{+}$, avec $B^{\pm}=\{x\in B:(x,c)=\pm\frac{1}{2}\}$.
On note $A$ la projection orthogonale de $B^{+}$ sur $V_{c}$. C'est
la $W_{c}$-orbite de $\frac{1}{n}(n-1,-1,\cdots,-1)\in V_{c}$. On
a $S_{n}^{W_{b}}=S_{b}^{W_{b}}[b]$ et $S_{n}^{W_{c}}=S_{c}^{W_{c}}[c]$.
D'après~\cite[§3]{Me88}, ou en adaptant la preuve de \cite[Chap. VI, §4, n°5 (IX)]{BoLie46}
en y remplaçant les polynômes symétriques élémentaires par des sommes
de puissances comme en \ref{subsec:GenAn}, on sait que
\[
\begin{array}{rclcrcl}
S_{n}^{W_{n}} & = & \mathbb{R}[\mathfrak{b}_{2},\mathfrak{b}_{4},\cdots,\mathfrak{b}_{2n}] & \text{avec} & \mathfrak{b}_{i} & = & \sum_{x\in B}x^{i}\\
S_{b}^{W_{b}} & = & \mathbb{R}[\mathfrak{b}_{2}^{\prime},\mathfrak{b}_{4}^{\prime},\cdots,\mathfrak{b}_{2(n-1)}^{\prime}] & \text{avec} & \mathfrak{b}_{i}^{\prime} & = & \sum_{y\in B^{\prime}}y^{i}\\
S_{c}^{W_{c}} & = & \mathbb{R}[\mathfrak{a}_{2},\mathfrak{a}_{3},\cdots,\mathfrak{a}_{n}] & \text{avec} & \mathfrak{a}_{i} & = & \sum_{z\in A}z^{i}
\end{array}
\]
Comme dans le cas de $A_{n}$, on vérifie que 
\[
(\beta-1)\mathfrak{r}_{2}^{B}=4(2n-1)b+2(2n-1)\quad\text{et}\quad(\gamma-1)\mathfrak{r}_{2}^{C}=8(n+1)c+(n^{2}+n)
\]
où $\mathfrak{r}_{2}^{B}=\sum_{\alpha\in R_{n}^{B}}\alpha^{2}$ et
$\mathfrak{r}_{2}^{C}=\sum_{\alpha\in R_{n}^{C}}\alpha^{2}$ sont
fixés par $W_{n}$, donc $b\in\boldsymbol{S}_{b}$ et $c\in\boldsymbol{S}_{c}$.
Puisque $B=\{\pm b\}\coprod B'$, $\mathfrak{b}_{2i}=2b^{2i}+\mathfrak{b}_{2i}^{\prime}$
et $\mathfrak{b}_{2i}^{\prime}\in\boldsymbol{S}_{b}$ pour tout $i$,
donc $\boldsymbol{S}_{b}=S_{n}^{W_{b}}$. D'autre part, $\mathfrak{b}_{2i}=2\mathfrak{b}_{2i}^{+}$
où $\mathfrak{b}_{i}^{+}=\sum_{y\in B^{+}}y^{i}$, donc $\mathfrak{b}_{2i}^{+}\in S_{n}^{W_{n}}$
pour tout $i$. Or $(\gamma-1)\mathfrak{b}_{2i}^{+}=\sum_{j=0}^{2i-1}\left(\begin{smallmatrix}2i\\
j
\end{smallmatrix}\right)(\frac{1}{2})^{2i-j}\mathfrak{b}_{j}^{+}$, donc $\mathfrak{b}_{i}^{+}\in\boldsymbol{S}_{c}$ pour tout $i$
par récurrence. Enfin $A=B^{+}-\frac{2}{n}c$, donc $\mathfrak{a}_{i}=\sum_{j=0}^{i}\left(\begin{smallmatrix}i\\
j
\end{smallmatrix}\right)\mathfrak{b}_{j}^{+}(\frac{-2}{n}c)^{i-j}\in\boldsymbol{S}_{c}$ pour tout $i$, et $\boldsymbol{S}_{c}=S_{n}^{W_{c}}$. 

\subsection{~}

Soient maintenant $W_{b}^{\circ}$, $W_{c}^{\circ}$ et $W_{c'}^{\circ}$
les stabilisateurs de $b$, $c$ et $c'$ dans $W_{n}^{\circ}$. On
note $\boldsymbol{S}_{b}^{\circ}$, $\boldsymbol{S}_{c}^{\circ}$,
et $\boldsymbol{S}_{c'}^{\circ}$ les sous-algèbres de $S_{n}^{W_{b}^{\circ}}$,
$S_{n}^{W_{c}^{\circ}}$ et $S_{n}^{W_{c'}^{\circ}}$ engendrées par
$S_{n}^{W_{n}^{\circ}}$ et son image par, respectivement, $\beta=\tau(b)$,
$\gamma=\tau(c)$ et $\gamma'=\tau(c')$. On vérifie que $W_{c}^{\circ}=W_{c}$,
donc $S_{n}^{W_{c}}=\boldsymbol{S}_{c}\subset\boldsymbol{S}_{c}^{\circ}\subset S_{n}^{W_{c}^{\circ}}=S_{n}^{W_{c}}$,
i.e.~$\boldsymbol{S}_{c}^{\circ}=S_{n}^{W_{c}^{\circ}}$. Puisque
$c$ et $c'$ sont conjugués sous $W_{n}$, on a de même $\boldsymbol{S}_{c'}^{\circ}=S_{n}^{W_{c'}^{\circ}}$.
Il reste à montrer que $\boldsymbol{S}_{b}^{\circ}=S_{n}^{W_{b}^{\circ}}$.
D'après~\cite[§3]{Me88} ou \cite[Chap. VI, §4, n°8 (IX)]{BoLie46},
$S_{n}^{W_{n}^{\circ}}=S_{n}^{W_{n}}[\mathfrak{d}]$ avec $\mathfrak{d}=\prod_{x\in B^{+}}x$,
et de même $S_{b}^{W_{b}^{\circ}}=S_{b}^{W_{b}}[\mathfrak{d}']$ avec
$\mathfrak{d}'=\prod_{x\in B^{\prime+}}x$ où $B'^{+}=B^{\prime}\cap B^{+}$.
On voit donc que 
\[
S_{n}^{W_{b}^{\circ}}=S_{b}^{W_{b}^{\circ}}[b]=S_{b}^{W_{b}}[\mathfrak{d}^{\prime},b]=S_{n}^{W_{b}}[\mathfrak{d}^{\prime}].
\]
Puisque $S_{n}^{W_{b}}=\boldsymbol{S}_{b}\subset\boldsymbol{S}_{b}^{\circ}$,
il reste à vérifier que $\mathfrak{d}^{\prime}\in\boldsymbol{S}_{b}^{\circ}$.
Or $B^{+}=\{b\}\cup B^{\prime+}$, donc $\mathfrak{d}=b\mathfrak{d}^{\prime}$,
et $(\beta-1)(\mathfrak{d})=(\beta(b)-b)\mathfrak{d}^{\prime}=\mathfrak{d}^{\prime}$
est dans $\boldsymbol{S}_{b}^{\circ}$ puisque $\mathfrak{d}$ est
dans $S_{n}^{W_{n}^{\circ}}$.

\subsection{Remarque\label{subsec:RemarqueDn}}

Soient $C$ et $C'$ les $W_{n}^{\circ}$-orbites de $c$ et $c'$.
Alors $\mathfrak{c}_{i}=\sum_{y\in C}y^{i}$ et $\mathfrak{c}_{i}^{\prime}=\sum_{y\in C'}y^{i}$
sont fixés par $W_{n}^{\circ}$ et échangés par tout élément de $W_{n}\setminus W_{n}^{\circ}$.
D'après ce qui précède, il existe un unique couple $(\mathfrak{c}_{i}^{0},\mathfrak{c}_{i}^{1})$
de $S_{n}^{W_{n}}$ tel que $\mathfrak{c}_{i}=\mathfrak{c}_{i}^{0}+\mathfrak{c}_{i}^{1}\mathfrak{d}$
et $\mathfrak{c}_{i}^{\prime}=\mathfrak{c}_{i}^{0}-\mathfrak{c}_{i}^{1}\mathfrak{d}$.
Pour $i=n,$ un calcul dans l'algèbre de polynômes $S_{n}=\mathbb{R}[x;x\in B_{+}]$
montre que $\mathfrak{c}_{n}^{1}=\frac{n!}{2}\neq0$, donc $\mathfrak{d}=\frac{1}{n!}(\mathfrak{c}_{n}-\mathfrak{c}_{n}^{\prime})$
et $S_{n}^{W_{n}^{\circ}}=S_{n}^{W_{n}}[\mathfrak{c}_{n}]=S_{n}^{W_{n}}[\mathfrak{c}_{n}^{\prime}]$.
Lorsque $n$ est impair, $C'=-C$, donc $\mathfrak{c}_{i}^{\prime}=(-1)^{i}\mathfrak{c}_{i}$
pour tout $i$, et $\mathfrak{c}_{n}=\frac{n!}{2}\mathfrak{d}=-\mathfrak{c}_{n}^{\prime}$. 

\section{Les cas $E_{6}$ et $E_{7}$}

Références: \cite[Chap. VI, \S 4, n°10, 11, 12]{BoLie46} 

\subsection{~}

On munit $V_{8}=\mathbb{R}^{8}$ du produit scalaire $((x_{i}),(y_{i}))=\sum x_{i}y_{i}$
et de la forme quadratique $q(x)=\frac{1}{2}(x,x)$. On note $I_{8}=\mathbb{Z}^{n}$,
$D_{8}=\{(x_{i})\in I_{8}:\sum x_{i}\equiv0\bmod2\}$, et $E_{8}=D_{8}+\mathbb{Z}(\frac{1}{2},\cdots,\frac{1}{2})=\{(x_{i})\in\mathbb{Z}^{8}\cup(\frac{1}{2}+\mathbb{Z})^{8}:\sum x_{i}\in2\mathbb{Z}\}$.
C'est un réseau unimodulaire pair, engendré par ses racines $R_{8}=\{\alpha\in E_{8}:q(\alpha)=1\}$,
de base 
\[
\begin{array}{ccccccccc}
\alpha_{1} & 1 & -1\\
\alpha_{2} &  & 1 & -1\\
\alpha_{3} &  &  & 1 & -1\\
\alpha_{4} &  &  &  & 1 & -1\\
\alpha_{5} &  &  &  &  & 1 & -1\\
\alpha_{6} &  &  &  &  &  & 1 & -1\\
\alpha_{7} &  &  &  &  &  & 1 & 1\\
\alpha_{8} & \frac{-1}{2} & \frac{-1}{2} & \frac{-1}{2} & \frac{-1}{2} & \frac{-1}{2} & \frac{-1}{2} & \frac{-1}{2} & \frac{-1}{2}
\end{array}
\]
et de diagramme de Dynkin 
\[
\xymatrix{\alpha_{1}\ar@{-}[r] & \alpha_{2}\ar@{-}[r] & \alpha_{3}\ar@{-}[r] & \alpha_{4}\ar@{-}[r] & \alpha_{5}\ar@{-}[r] & \alpha_{7}\ar@{-}[r] & \alpha_{8}\\
 &  &  &  & \alpha_{6}\ar@{-}[u]
}
\]

\subsection{~}

On note $V_{7}=\{(x_{i}):\sum x_{i}=0\}$ l'orthogonal de $\alpha_{8}$
dans $V_{8}$. Alors $E_{7}=V_{7}\cap E_{8}$ est un réseau pair d'indice
$2$ dans son dual $E_{7}^{\ast}\subset V_{7}$, qui est engendré
par ses racines $R_{7}=V_{7}\cap R_{8}=\{\alpha\in E_{7}:q(\alpha)=1\}$,
de base et diagramme de Dynkin\xyR{.8pc} \xyC{.8pc}
\[
\begin{array}{ccccccccc}
\alpha_{1} & 1 & -1\\
\alpha_{2} &  & 1 & -1\\
\alpha_{3} &  &  & 1 & -1\\
\alpha_{4} &  &  &  & 1 & -1\\
\alpha_{5} &  &  &  &  & 1 & -1\\
\alpha_{6} &  &  &  &  &  & 1 & -1\\
\alpha_{7}^{\prime} & \frac{-1}{2} & \frac{-1}{2} & \frac{-1}{2} & \frac{-1}{2} & \frac{1}{2} & \frac{1}{2} & \frac{1}{2} & \frac{1}{2}
\end{array}\quad\xymatrix{\alpha_{1}\ar@{-}[r] & \alpha_{2}\ar@{-}[r] & \alpha_{3}\ar@{-}[r] & \alpha_{4}\ar@{-}[r] & \alpha_{5}\ar@{-}[r] & \alpha_{6}\\
 &  &  & \alpha_{7}^{\prime}\ar@{-}[u]
}
\]
La plus haute racine est $\alpha_{1}+2\alpha_{2}+3\alpha_{3}+4\alpha_{4}+3\alpha_{5}+2\alpha_{6}+2\alpha'_{7}$,
et le vecteur 
\[
a=\frac{1}{2}\left(\frac{3}{2},\frac{-1}{2},\cdots,\frac{-1}{2},\frac{3}{2}\right)\in E_{7}^{\ast}
\]
est l'unique élément dominant minuscule de $E_{7}^{\ast}$, orthogonal
à toutes les racines simples sauf $\alpha_{1}$. On a $q(a)=\frac{3}{4}$.
On note $W_{7}$ le groupe de Weyl de $R_{7}$, de cardinal $\left|W_{7}\right|=2^{10}\cdot3^{4}\cdot5\cdot7$;
c'est aussi le groupe orthogonal $O(E_{7})$ de $E_{7}$. 

On note $A$ la $W_{7}$-orbite de $a$ dans $V_{7}$. 

\subsection{~}

On note $V_{6}$ l'orthogonal de $a$ dans $V_{7}$. Alors $E_{6}=V_{6}\cap E_{7}$
est un réseau pair d'indice $3$ dans son dual $E_{6}^{\ast}\subset V_{6}$,
qui est engendré par ses racines $R_{6}=V_{6}\cap R_{7}=\{\alpha\in E_{6}:q(\alpha)=1\}$,
de base et diagramme de Dynkin
\[
\begin{array}{ccccccccc}
\alpha_{2} &  & 1 & -1\\
\alpha_{3} &  &  & 1 & -1\\
\alpha_{4} &  &  &  & 1 & -1\\
\alpha_{5} &  &  &  &  & 1 & -1\\
\alpha_{6} &  &  &  &  &  & 1 & -1\\
\alpha_{7}^{\prime} & \frac{-1}{2} & \frac{-1}{2} & \frac{-1}{2} & \frac{-1}{2} & \frac{1}{2} & \frac{1}{2} & \frac{1}{2} & \frac{1}{2}
\end{array}\quad\xymatrix{\alpha_{2}\ar@{-}[r] & \alpha_{3}\ar@{-}[r] & \alpha_{4}\ar@{-}[r] & \alpha_{5}\ar@{-}[r] & \alpha_{6}\\
 &  & \alpha_{7}^{\prime}\ar@{-}[u]
}
\]
La plus haute racine est $\alpha_{2}+2\alpha_{3}+3\alpha_{4}+2\alpha_{7}^{\prime}+2\alpha_{5}+\alpha_{6}$,
et les vecteurs
\begin{align*}
b^{+} & =\frac{1}{3}\left(\frac{-3}{2},\frac{5}{2},\frac{-1}{2},\cdots,\frac{-1}{2},\frac{3}{2}\right)\\
b^{-} & =\frac{1}{3}\left(\frac{-3}{2},\frac{1}{2},\cdots,\frac{1}{2},\frac{-5}{2},\frac{3}{2}\right)
\end{align*}
sont les éléments dominants minuscules de $E_{6}^{\ast}$, orthogonaux
à toutes les racines simples sauf respectivement $\alpha_{2}$ et
$\alpha_{6}$. On a $q(b^{\pm})=\frac{2}{3}$. On note $W_{6}$ le
groupe de Weyl de $R_{6}$. C'est le stabilisateur de $a$ dans $W_{7}$,
il est de cardinal $\left|W_{6}\right|=2^{7}\cdot3^{4}\cdot5$, donc
$\left|A\right|=2^{3}\cdot7=56$. Le groupe orthogonal de $E_{6}$
est $O(E_{6})=W_{6}\times\{\pm1\}$.

On note $B^{+}$ et $B^{-}$ les $W_{6}$-orbites de $b^{+}$ et $b^{-}$
dans $V_{6}$. 

\subsection{~}

Pour $\epsilon\in\{+,-\}$, on note $V_{5}^{\epsilon}$ l'orthogonal
de $b^{\epsilon}$ dans $V_{6}$. Alors $E_{5}^{\epsilon}=V_{5}^{\epsilon}\cap E_{6}$
est un réseau pair d'indice $4$ dans son dual $E_{5}^{\epsilon\ast}\subset V_{5}^{\epsilon}$,
qui est engendré par ses racines $R_{5}^{\epsilon}=V_{5}^{\epsilon}\cap R_{6}=\{\alpha\in E_{5}^{\epsilon}:q(\alpha)=1\}$.
Pour $\epsilon=+$, une base de $R_{5}^{+}$ est 
\[
\begin{array}{ccccccccc}
\alpha_{3} &  &  & 1 & -1\\
\alpha_{4} &  &  &  & 1 & -1\\
\alpha_{5} &  &  &  &  & 1 & -1\\
\alpha_{6} &  &  &  &  &  & 1 & -1\\
\alpha_{7}^{\prime} & \frac{-1}{2} & \frac{-1}{2} & \frac{-1}{2} & \frac{-1}{2} & \frac{1}{2} & \frac{1}{2} & \frac{1}{2} & \frac{1}{2}
\end{array}\quad\xymatrix{\alpha_{3}\ar@{-}[r] & \alpha_{4}\ar@{-}[r] & \alpha_{5}\ar@{-}[r] & \alpha_{6}\\
 & \alpha_{7}^{\prime}\ar@{-}[u]
}
\]
la plus haute racine est $\alpha_{3}+\alpha_{7}^{\prime}+\alpha_{6}+2\alpha_{4}+2\alpha_{5}$
et les vecteurs 
\[
c^{+}={\textstyle \frac{1}{4}\left(-1,-1,1,\cdots,1,-3,1\right)\quad\text{et}\quad d^{+}=\frac{1}{8}\left(-3,-3,7,-1,\cdots,-1,3\right)}
\]
de $E_{5}^{+\ast}$ sont dominants minuscules, orthogonaux à toutes
les racines simples sauf respectivement $\alpha_{6}$ et $\alpha_{3}$.
Pour $\epsilon=-$, une base de $R_{5}^{-}$ est 
\[
\begin{array}{ccccccccc}
\alpha_{2} &  & 1 & -1\\
\alpha_{3} &  &  & 1 & -1\\
\alpha_{4} &  &  &  & 1 & -1\\
\alpha_{5} &  &  &  &  & 1 & -1\\
\alpha_{7}^{\prime} & \frac{-1}{2} & \frac{-1}{2} & \frac{-1}{2} & \frac{-1}{2} & \frac{1}{2} & \frac{1}{2} & \frac{1}{2} & \frac{1}{2}
\end{array}\quad\xymatrix{\alpha_{2}\ar@{-}[r] & \alpha_{3}\ar@{-}[r] & \alpha_{4}\ar@{-}[r] & \alpha_{5}\\
 &  & \alpha_{7}^{\prime}\ar@{-}[u]
}
\]
la plus haute racine est $\alpha_{2}+\alpha_{5}+\alpha_{7}^{\prime}+2\alpha_{3}+2\alpha_{4}$,
et les vecteurs
\[
c^{-}={\textstyle \frac{1}{4}\left(-1,3,-1,\cdots,-1,1,1\right)\quad\text{et}\quad d^{-}=\frac{1}{8}\left(-3,1,\cdots,1,-7,3,3\right)}
\]
de $E_{5}^{-\ast}$ sont dominants minuscules, orthogonaux à toutes
les racines simples sauf respectivement $\alpha_{2}$ et $\alpha_{5}$.
Pour $\epsilon\in\{+,-\}$, on note $W_{5}^{\epsilon}$ le groupe
de Weyl de $R_{5}^{\epsilon}$. C'est le stabilisateur de $b^{\epsilon}$
dans $W_{6}$, il est de cardinal $\left|W_{5}^{\epsilon}\right|=2^{7}\cdot3\cdot5$,
donc $\left|B^{\epsilon}\right|=2^{3}=27$. Le groupe orthogonal de
$E_{5}^{\epsilon}$ est $O(E_{5}^{\epsilon})=W_{5}^{\epsilon}\times\{\pm1\}$. 

On note $C^{\epsilon}$ et $D^{\epsilon}$ les $W_{5}^{\epsilon}$-orbites
de $c^{\epsilon}$ et $d^{\epsilon}$ dans $V_{5}^{\epsilon}$. Le
stabilisateur de $c^{\epsilon}$ (resp.~$d^{\epsilon}$) dans $W_{5}^{\epsilon}$
est le groupe de Weyl d'un système de racines de type $D_{4}$ (resp.~$A_{4}$),
il est de cardinal $2^{6}\cdot3$ (resp.~$2^{3}\cdot3\cdot5$), donc
$\left|C^{\epsilon}\right|=10$ et $\left|D^{\epsilon}\right|=16$. 

\subsection{~}

Revenons à $E_{7}$ et à la $W_{7}$-orbite $A$ de $a$. Puisque
$\left|A\right|=56$ et $W_{7}$ contient $\mathfrak{S}_{8}$ et $\{\pm1\}$,
$A$ est l'ensemble de tous les vecteurs obtenus par permutation des
coordonnées de $\pm\frac{1}{2}(\frac{3}{2},\frac{3}{2},\frac{-1}{2},\cdots,\frac{-1}{2})$.
Considérons sur $A\times A$ la fonction $W_{7}$-invariante définie
par $d(x,x')=q(x-x')=\frac{3}{2}-(x,x')$. On vérifie qu'elle est
à valeur dans $\{0,1,2,3\}$, que c'est une distance, et que $d(x,x')+d(x,-x')=3$.
En fait, l'ensemble $A$ muni de la relation d'incidence $x\sim x'\iff d(x,x')=1$
est le graphe considéré dans \cite{Co90}\footnote{Avec les notations de \cite{Co90}, l'identification de $\Omega=P\cup P^{\ast}$
avec $A$ se fait en associant à $p_{\alpha}\in P$ (resp.~$p_{\alpha}^{\ast}\in P^{\ast}$)
le vecteur de coordonnées $x_{i}=\frac{3}{4}$ (resp.~$-\frac{3}{4}$)
pour $i\in\alpha$ et $x_{i}=-\frac{1}{4}$ (resp.~$\frac{1}{4}$)
pour $i\notin\alpha$, où $\alpha$ décrit l'ensemble $\Phi^{(2)}$
des parties à deux élements de $\Phi=\{1,\cdots,8\}$.}, dont $d$ est la distance, et dont $W_{7}$ est le groupe des automorphismes. 

\subsection{~}

On en déduit d'abord une partition de $A$ en quatre parties stables
sous $W_{6}$, échangées deux à deux par $-1\in W_{7}$: les points
fixes $\{a\}$ et $\{-a\}$, et deux parties à $27$ éléments, $A^{+}$
et $A^{-}=-A^{+}$, qui sont respectivement définies par
\[
A^{+}=\{x\in A:d(a,x)=1\}\quad\text{et}\quad A^{-}=\{x\in A:d(a,x)=2\}.
\]
D'autre part, on vérifie que 
\[
\begin{array}{rccccl}
a^{+} & := & b^{+}+{\textstyle \frac{1}{3}a} & = & \frac{1}{4}\left(-1,3,-1,\cdots,-1,3\right) & \in A^{+}\\
\text{et}\quad a^{-} & := & b^{-}-{\textstyle \frac{1}{3}a} & = & \frac{1}{4}\left(-3,1,\cdots,1,-3,1\right) & \in A^{-}.
\end{array}
\]
Pour $\epsilon\in\{+,-\}$ et puisque $\left|B^{\epsilon}\right|=27$,
il en résulte que $A^{\epsilon}=B^{\epsilon}+\frac{\epsilon}{3}a$
est la $W_{6}$-orbite de $a^{\epsilon}$, et que $B^{\epsilon}$
est la projection orthogonale de $A^{\epsilon}\subset V_{7}$ sur
$V_{6}$. En particulier, les $W_{6}$-orbites $B^{+}$ et $B^{-}$
sont échangées par $-1\in O(E_{6})$. Puisque la restriction de $d$
à $A^{\epsilon}$ est à valeur dans $\{0,1,2\}$, on obtient ensuite
une décomposition de $A^{\epsilon}$ en trois parties stables sous
$W_{5}^{\epsilon}$: le point fixe $\{a^{\epsilon}\}$, et les sous-ensembles
\[
\begin{array}{rclcrcl}
C_{7}^{\epsilon} & = & \left\{ x\in A^{\epsilon}:d(a^{\epsilon},x)=2\right\}  & \quad\text{et}\quad & D_{7}^{\epsilon} & = & \left\{ x\in A^{\epsilon}:d(a^{\epsilon},x)=1\right\} \end{array}
\]
qui ont respectivement $10$ et $16$ éléments. Leur projection orthogonale
sur $V_{6}$ donne une décomposition analogue de $B^{\epsilon}$ en
trois parties stables sous $W_{5}^{\epsilon}$: le point fixe $\{b^{\epsilon}\}$
et des sous-ensembles $C_{6}^{\epsilon}=C_{7}^{\epsilon}-\frac{\epsilon}{3}a$
et $D_{6}^{\epsilon}=D_{7}^{\epsilon}-\frac{\epsilon}{3}a$, caractérisés
par
\[
\begin{array}{rclcrcl}
C_{6}^{\epsilon} & = & \left\{ y\in B^{\epsilon}:(y,b^{\epsilon})={\textstyle \frac{-2}{3}}\right\}  & \quad\text{et}\quad & D_{6}^{\epsilon} & = & \left\{ y\in B^{\epsilon}:(y,b^{\epsilon})={\textstyle \frac{1}{3}}\right\} \end{array}.
\]
On vérifie enfin que 
\[
\begin{array}{rclcrcl}
c^{\epsilon}-{\textstyle \frac{1}{2}b^{\epsilon}+\frac{\epsilon}{3}a} & \in & C_{7}^{\epsilon} & \quad\text{et}\quad & d^{\epsilon}+\frac{1}{4}b^{\epsilon}+\frac{1}{3}\epsilon a & \in & D_{7}^{\epsilon},\end{array}
\]
donc
\[
\begin{array}{rclcrcl}
c^{\epsilon}-{\textstyle \frac{1}{2}b^{\epsilon}} & \in & C_{6}^{\epsilon} & \quad\text{et}\quad & d^{\epsilon}+\frac{1}{4}b^{\epsilon} & \in & D_{6}^{\epsilon}.\end{array}
\]
Puisque $\left|C^{\epsilon}\right|=10$ et $\left|D^{\epsilon}\right|=16$,
il en résulte que $C_{6}^{\epsilon}=C^{\epsilon}-\frac{1}{2}b^{\epsilon}$
et $D_{6}^{\epsilon}=D^{\epsilon}+\frac{1}{4}b^{\epsilon}$. On en
déduit que les $W_{5}^{\epsilon}$-orbites $C^{\epsilon}$ et $D^{\epsilon}$
de $V_{5}^{\epsilon}$ sont les projections orthogonales des $W_{5}^{\epsilon}$-orbites
que sont finalement $C_{7}^{\epsilon}$ et $D_{7}^{\epsilon}$ dans
$V_{7}$, ou $C_{6}^{\epsilon}$ et $D_{6}^{\epsilon}$ dans $V_{6}$.

\subsection{~}

Soient $S_{5}^{\epsilon}\subset S_{6}\subset S_{7}$ les algèbres
symétriques de $V_{5}^{\epsilon}\subset V_{6}\subset V_{7}$. On a
donc $S_{7}^{W_{6}}=S_{6}^{W_{6}}[a]$ et $S_{6}^{W_{5}^{\epsilon}}=(S_{5}^{\epsilon})^{W_{5}^{\epsilon}}[b^{\epsilon}]$.
Il résulte de \cite{Me88}\footnote{Avec les notations de \cite[§ 2.1]{Me88} pour $E_{6}$, $(y_{1},x_{1},\cdots,x_{6},y_{2})$
est la base duale de la base canonique de $V_{8}$, $-\frac{1}{2}(y_{1}-y_{2})+\frac{1}{6}S_{1}-x_{6}$
est la forme linéaire $(b^{-},-)$ sur $V_{8}$, et les générateurs
$I_{k}$ de Mehta pour $S_{6}^{W_{6}}$ sont donc nos $\mathfrak{b}_{k}^{-}=(-1)^{k}\mathfrak{b}_{k}^{+}$
pour $k\in\{2,5,6,8,9,12\}$.} et de la remarque~\ref{subsec:RemarqueDn} que $S_{6}^{W_{6}}$
et $(S_{5}^{\epsilon})^{W_{5}^{\epsilon}}$ sont les algèbres graduées
respectivement engendrées par 
\[
\mathfrak{b}_{i}^{\epsilon}=\sum_{b\in B^{\epsilon}}b^{i}\quad\text{pour }i\in\{2,5,6,8,9,12\},
\]
\[
\mathfrak{c}_{i}^{\epsilon}=\sum_{c\in C^{\epsilon}}c^{i}\quad\text{pour }i\in\{2,4,6,8\}\quad\text{et}\quad\mathfrak{d}_{i}^{\epsilon}=\sum_{d\in D^{\epsilon}}d^{i}\quad\text{pour }i=5.
\]
Soient $\tau=\tau(a)$ et $\tau^{\epsilon}=\tau(b^{\epsilon})$. Notons
$\boldsymbol{S}_{7}$ la sous-algèbre de $S_{7}^{W_{6}}$ engendrée
par $S_{7}^{W_{7}}$ et $\tau(S_{7}^{W_{7}})$, et $\boldsymbol{S}_{6}^{\epsilon}$
la sous-algèbre de $S_{6}^{W_{5}^{\epsilon}}$ engendrée par $S_{6}^{W_{6}}$
et $\tau^{\epsilon}(S_{6}^{W_{6}})$. On se propose de montrer que
$\boldsymbol{S}_{7}=S_{7}^{W_{6}}$ et $\boldsymbol{S}_{6}^{\epsilon}=S_{6}^{W_{5}^{\epsilon}}$.
Il suffit donc de montrer que $(1)$ $a$ et tous les $\mathfrak{b}_{i}^{+}$
sont dans $\boldsymbol{S}_{7}$, et $(2)$ $b^{\epsilon}$ et tous
les $\mathfrak{c}_{i}^{\epsilon}$, $\mathfrak{d}_{i}^{\epsilon}$
sont dans $\boldsymbol{S}_{6}^{\epsilon}$. 

\subsection{~}

Posons $\mathfrak{r}_{7,2}=\sum_{\alpha\in R_{7}}\alpha^{2}$ et $\mathfrak{r}_{6,2}=\sum_{\alpha\in R_{6}}\alpha^{2}$,
qui sont respectivement dans $S_{7}^{W_{7}}$ et $S_{6}^{W_{6}}$.
On vérifie comme dans les cas précédents les relations
\[
(\tau-1)(\mathfrak{r}_{7,2})=54+72a\quad\text{et}\quad(\tau^{\epsilon}-1)(\mathfrak{r}_{6,2})=32+48b^{\epsilon},
\]
donc $a\in\boldsymbol{S}_{7}$ et $b^{\epsilon}\in\boldsymbol{S}_{6}^{\epsilon}$. 

\subsection{~}

Notons $\mathfrak{a}_{i}=\sum_{x\in A}x^{i}$, un élément de degré
$i$ dans $S_{7}^{W_{7}}$. Comme 
\[
A=\{a\}{\textstyle \coprod}A^{+}{\textstyle \coprod}(-A^{+}){\textstyle \coprod}\{-a\}
\]
on voit que $\mathfrak{a}_{i}=0$ pour $i\equiv1\bmod2$, tandis que
$\mathfrak{a}_{i}=2(a^{i}+\mathfrak{a}_{i}^{+})$ pour $i\equiv0\bmod2$,
avec $\mathfrak{a}_{i}^{+}=\sum_{x\in A^{+}}x^{i}$. Comme $a^{i}$
et $(\tau-1)(a^{i})=(a+\frac{3}{2})^{i}-a^{i}$ sont dans $\boldsymbol{S}_{7}$,
on en déduit que $\mathfrak{a}_{2i}^{+}$ et $(\tau-1)(\mathfrak{a}_{2i}^{+})$
sont dans $\boldsymbol{S}_{7}$ pour tout $i\in\mathbb{N}$. Or $(a,A^{+})=\frac{1}{2}$,
donc
\[
(\tau-1)\mathfrak{a}_{2i}^{+}=\sum_{x\in A^{+}}(x+{\textstyle \frac{1}{2}})^{2i}-x^{2i}=\sum_{j=0}^{2i-1}\left(\begin{smallmatrix}2i\\
j
\end{smallmatrix}\right)2^{j-2i}\mathfrak{a}_{j}^{+}
\]
et $\mathfrak{a}_{i}^{+}\in\boldsymbol{S}_{7}$ pour tout $i\in\mathbb{N}$
par récurrence. Enfin, $B^{+}=A^{+}-\frac{1}{3}a$, donc
\[
\mathfrak{b}_{i}^{+}=\sum_{x\in A^{+}}(x-{\textstyle \frac{1}{3}}a)^{i}=\sum_{j=0}^{i}\left(\begin{smallmatrix}i\\
j
\end{smallmatrix}\right)(-{\textstyle \frac{1}{3}a)^{i-j}\mathfrak{a}_{j}^{+}\in\boldsymbol{S}_{7}.}
\]
On a donc bien $\boldsymbol{S}_{7}=S_{7}^{W_{6}}$. 

\subsection{~}

Comme $B^{\epsilon}=\{b^{\epsilon}\}{\textstyle \coprod}C_{6}^{\epsilon}{\textstyle \coprod}D_{6}^{\epsilon}$,
on voit que $\mathfrak{\mathfrak{b}}_{i}^{\epsilon}=(b^{\epsilon})^{i}+\mathfrak{C}_{i}^{\epsilon}+\mathfrak{D}_{i}^{\epsilon}$
où
\[
\mathfrak{C}_{i}^{\epsilon}=\sum_{y\in C_{6}^{\epsilon}}y^{i}\quad\text{et}\quad\mathfrak{D}_{i}^{\epsilon}=\sum_{z\in D_{6}^{\epsilon}}z^{i}.
\]
Comme $(b^{\epsilon})^{i}$ et $(\tau^{\epsilon}-1)(b^{\epsilon})^{i}=(b^{\epsilon}+\frac{4}{3})^{i}-(b^{\epsilon})^{i}$
sont dans $\boldsymbol{S}_{6}^{\epsilon}$, on en déduit que $\mathfrak{C}_{i}^{\epsilon}+\mathfrak{D}_{i}^{\epsilon}$
et $(\tau^{\epsilon}-1)(\mathfrak{C}_{i}^{\epsilon}+\mathfrak{D}_{i}^{\epsilon})$
sont dans $\boldsymbol{S}_{6}^{\epsilon}$. Or $(b^{\epsilon},C_{6}^{\epsilon})=\frac{-2}{3}$
et $(b^{\epsilon},D_{6}^{\epsilon})=\frac{1}{3}$, donc 
\[
(\tau^{\epsilon}-1)(\mathfrak{C}_{i}^{\epsilon}+\mathfrak{D}_{i}^{\epsilon})=\sum_{j=0}^{i-1}\left(\begin{smallmatrix}i\\
j
\end{smallmatrix}\right)\cdot\left(({\textstyle -\frac{2}{3}})^{i-j}\mathfrak{C}_{j}^{\epsilon}+({\textstyle \frac{1}{3}})^{i-j}\mathfrak{D}_{j}^{\epsilon}\right).
\]
La matrice $\left(\begin{smallmatrix}1 & 1\\
\frac{-2}{3} & \frac{1}{3}
\end{smallmatrix}\right)$ étant inversible, on en déduit par récurrence sur $i$ que $\mathfrak{C}_{i}^{\epsilon}$
et $\mathfrak{D}_{i}^{\epsilon}$ sont dans $\boldsymbol{S}_{6}^{\epsilon}$
pour tout $i\in\mathbb{N}$. Comme enfin $C^{\epsilon}=C_{6}^{\epsilon}+\frac{1}{2}b^{\epsilon}$
et $D^{\epsilon}=D_{6}^{\epsilon}-\frac{1}{4}b^{\epsilon}$, 
\begin{align*}
\mathfrak{c}_{i}^{\epsilon} & =\sum_{y\in C_{6}^{\epsilon}}(y+{\textstyle \frac{1}{2}}b^{\epsilon})^{i}=\sum_{j=0}^{i}\left(\begin{smallmatrix}i\\
j
\end{smallmatrix}\right)({\textstyle \frac{1}{2}b^{\epsilon})^{i-j}\mathfrak{C}_{j}^{\epsilon}}\\
\text{et}\quad\mathfrak{d}_{i}^{\epsilon} & =\sum_{z\in D_{6}^{\epsilon}}(z-{\textstyle \frac{1}{4}}b^{\epsilon})^{i}=\sum_{j=0}^{i}\left(\begin{smallmatrix}i\\
j
\end{smallmatrix}\right)({\textstyle \frac{-1}{4}b^{\epsilon})^{i-j}\mathfrak{D}_{j}^{\epsilon}}
\end{align*}
sont également dans $\boldsymbol{S}_{6}^{\epsilon}$: on a donc bien
$\boldsymbol{S}_{6}^{\epsilon}=S_{6}^{W_{5}^{\epsilon}}$. 

\bibliographystyle{plain}
\bibliography{/home/christophe/Maths/MyBib}

\end{document}